\documentstyle[12pt]{article}
\newcommand\CC{\hbox{C\kern -.58em {\raise .54ex \hbox{$\scriptscriptstyle |$}}
  \kern-.55em {\raise .53ex \hbox{$\scriptscriptstyle |$}} }}
\newcommand\qed{\hfill$\sqcap\kern-8.0pt\hbox{$\sqcup$}$}
\newcommand\NN{\hbox{I\kern-.2em\hbox{N}}}
\newcommand\RR{\hbox{I\kern-.2em\hbox{R}}}
\newcommand\sRR{{\sl \hbox{I\kern-.2em\hbox{R}}}}
\newcommand\QQ{\hbox{I\kern-.53em\hbox{Q}}}
\newcommand\ZZ{{{\rm Z}\kern-.28em{\rm Z}}}

\newtheorem{thm}{Theorem}[section]

\newtheorem{lem}[thm]{Lemma}
\newtheorem{rem}[thm]{Remark}

\newtheorem{df}[thm]{Definition}
\newtheorem{exm}[thm]{Example}

\def\fnote#1{\footnote}

\def\abstract{\if@twocolumn
\section*{Abstract}
\else \small {\bf ABSTRACT\vspace{-.5em}\vspace{0pt}} \quotation
\fi}

\begin{document}
\title{{\LARGE{A constructive method for linear extensions of Zadeh's fuzzy order}}}
\author{\small{Abdelkader STOUTI}}
\date{\small{Laboratory of Mathematics and Applications, Faculty of Sciences and Techniques,
University Sultan Moulay Slimane, P.O. Box 523, 23000 Beni-Mellal, MOROCCO. \\E-mail: stout@fstbm.ac.ma}}\maketitle

\bigskip

\begin{center}
\begin{abstract}
\noindent  In this paper, we give a constructive method for linear
extensions of Zadeh's fuzzy orders. We also characterize Zadeh's
fuzzy orders by their linear extensions.
\end{abstract}
\end{center}

\bigskip

\medskip

\medskip


\noindent \baselineskip=20pt

\section{Introduction and preliminaries}
In 1971, Zadeh \cite{Zadeh71} introduced the notion of fuzzy
order.  Since then, the theory of fuzzy order was developed by
many authors (for example see: [1-4] and [6, 8, 9, 12, 13]). In
this paper, we first establish a constructive method for linear
extensions of Zadeh's fuzzy orders. Secondly, we prove that every
Zadeh's fuzzy order defined on a finite set is the fuzzy
intersection of its all linear Zadeh's fuzzy extensions.

\bigskip

This paper is organized as follows. In the first section, we
recall some well know definitions and results. In the seond
section, we first give the key result of the present paper (see
Lemma 2.2). For the case of finite sets, a constructive method for
linear extensions of Zadeh's fuzzy orders is given in the second
section (see Theorem 2.1). In the third section, we characterize
Zadeh's fuzzy orders by their linear extensions (see Theorem 3.1).
In the forth section, we give some examples of building linear
extensions of Zadeh's fuzzy orders.

\medskip

Next, we will recall some well know definitions and results.

\medskip

Let $X$ be a nonempty set. A fuzzy subset $A$ of $X$ is
characterized by its membership function $A : X \rightarrow [0,
1]$ and $A(x)$ is interpreted as the degree of membership of the
element $x$ in the fuzzy subset $A$ for each $x\in X.$

\begin{df} \cite{Zadeh71}. Let $X$ be a nonempty set.
A fuzzy relation $r$  on  $X$ is a function $r :X\times X
\longrightarrow [0 , 1]$. For every $x, y\in X,$ the value
$r(x,y)$ is called the grade of membership of $(x,y)$ in $r$ and
means how far x and y are related under $r.$
\end{df}

In \cite{Zadeh71}, Zadeh gave the following definition of fuzzy
order.

\begin{df} \cite{Zadeh71}. Let $X$ be a nonempty crisp set. A Zadeh's fuzzy
order on $X$ is a fuzzy subset $r$ of $X\times X$ satisfying the
following three properties:

(i) for all $x \in X,$ $r(x, x)=1$, (Z-fuzzy reflexivity);

(ii) for all $x, y \in X,$ $x\not=y$ and $r(x, y) > 0 $ imply
$r(y, x)=0$ (Z-fuzzy antisymmetry);

iii) for all $x, z \in X,$  $r(x, z) \geq \sup_{y \in X}\left[
\min\{ r(x, y), r(y, z)\} \right],$ (Z-fuzzy transitivity).
\end{df}

 A nonempty set $X$ with a Zadeh's fuzzy order $r$
defined on it is called Zadeh's--fuzzy ordered set (for short,
Z-foset) and we denote it by $(X, r).$

If $Y$ is a subset of an Z-foset $(X, r)$, then the Z-fuzzy order
$r$ is a Z-fuzzy order on $Y$ and is called the induced fuzzy
order.

A Z-fuzzy order $r$ is linear (total) on $X$ if for every $x, y
\in X,$ we have $r(x, y) > 0$ or $r(y, x) > 0.$ A  Z-fuzzy ordered
set $(X, r)$ in which $r$ is linear is called a Z-fuzzy chain. If
 there is at least two elements $x,y \in X$ such that $r(x,y)=r(y, x)=0,$
then the Z-fuzzy order $r$ is said to be a partial Z-fuzzy order.

\bigskip

Let $(X, r)$ be a Z-fuzzy ordered set and $A$ be a subset of $X$.

(a) An element $u \in X$ is a $r$-upper bound of $A$ if $r(x, u) >
0$ for all $x \in A.$ If $u$ is $r$-upper bound of $A$ and $u \in
A,$ then $u$ is called a greatest element of $A.$ The $r$-lower
bound and least element are defined analogously.

(b) An element $m \in A$ is called a maximal element of $A$ if
$r(m, x) > 0$ for some $x \in A,$ then $x=m.$ Minimal elements are
defined similarly.

(c) An element $s \in X$ is the $r$-supremum of $A$ if $s$ is a
$r$-upper bound of $A$ and for all $r$-upper bound $u$ of $A,$ we
have $r(s, u) > 0.$ When $s$ exists, we shall write
$s=\sup_{r}(A).$ Similarly, $l \in X$ is the $r$-infimum of $A$ if
$l$ is a $r$-lower bound of $A$ and for all $r_{\alpha}$-lower
bound $k$ of $A,$ we have $r(k, l) > 0.$ When $l$ exists we shall
write $l=\inf_{r}(A).$

\begin{df}
Let $X$ be a nonempty set and $r$, $r'$ be two Zadeh's fuzzy
orders on $X$. We say that $r'$ is an extension of $r$ if
$r(x,y)\leq r'(x,y),$ for every $(x, y) \in X^{2}.$
\end{df}

\begin{df}
Let $(X, r)$ be a nonempty Zadeh's fuzzy ordered set and let $a,
b$ be two elements of $X.$ We say that $a$ and $b$ are
incomparable in $(X, r)$ if $r(a,b)=r(b,a)=0.$
\end{df}

\noindent Let $x, y \in \RR.$ Then, we set $\max\{x, y\}=x\vee y
\hbox{ and } \min\{x, y\}=x\wedge y.$

\begin{exm}

\bigskip

 \noindent 1.  Let $r_{1}$ and $r_{2}$ be the two fuzzy
relations defined on $\RR$ by:
 $$r_{1}(x,y) =
\left\{\begin{array}{c}
1 ~,~ if ~~ x=y ;~~~~~~~~~~~~~\\
\min(1 , \frac{y-x}{2}), ~~if  ~~ x<y ~~~ \\
0, ~~ ~~if  ~~ x> y~~~~~~~~~~~~~
\end{array}
\right.$$ and  $$r_{2}(x,y) = \left\{
\begin{array}{c}
1 \quad \hbox{ if } \quad x=y\\
0,75 \quad \hbox{ if } \quad x > y\\
0 \quad \hbox{ if } \quad x < y.
\end{array}
\right.$$
Then, $r_{1}$ and $r_{2}$ are two Zadeh's fuzzy orders
on $\RR.$

 \noindent 2. Let $X=\{a, b, c\}$ and $r$ be the fuzzy relation
 defined on $X$ by the following table:

\begin{center}

\begin{tabular}{|c|c|c|c|}
  \hline
& a &b & c \\
 \hline
a & 1 & 0 & $\gamma$ \\
b & 0 & 1 & 0 \\
c & 0 & 0 & 1 \\
  \hline
\end{tabular}
\end{center}
Then, $r$ is a  Zadeh's fuzzy order on $X.$
\medskip

\noindent 3. Let $X=\{a, b, c, d\}$ and $r$ be the fuzzy relation
 defined on $X$ by the following  matrix:

\begin{center}

\begin{tabular}{|c|c|c|c|c|}
  \hline
& a & b & c & d \\
 \hline
a & 1 & 0 & $\gamma$ & $\lambda$ \\
b & 0 & 1 & $\alpha$ & $\beta$ \\
c & 0 & 0 & 1 & 0 \\
d & 0 & 0 & 0 & 1 \\
  \hline
\end{tabular}
\end{center}
Then, $r$ is a  Zadeh's fuzzy order on $X.$
\end{exm}
\section{A constructive method for linear extensions of Zadeh's fuzzy orders defined on finite sets}
In this section,  we will show how one can constructs a linear
extension of a partial Zadeh's fuzzy order defined on a finite
nonempty set $X.$ More precisely, we will prove the following
result.

\begin{thm}
Every Zadeh's fuzzy order on finite nonempty set $X$ can be
extended to a linear Zadeh's fuzzy order on $X$.
\end{thm}

In order to prove Theorem 2.1, we will need the following
technical lemma.

\begin{lem}  Let $(X, r)$ be a nonempty Zadeh's
 fuzzy ordered set and let $a, b$ be two elements of $X$ such that
 $r(b, a)=0.$ Then, there exists at least a Zadeh's fuzzy order $r'$ on $X$
which extends $r$ such that $r'(a,b)=1$ and $r'(b,a)=0.$
\end{lem}

\noindent{\bf Proof}. Let $(X, r)$ be a Zadeh's fuzzy ordered set
and let $a, b\in X$ such that $r(b,a)=0.$ Let $s$ be the fuzzy
relation defined on $X$ by setting
 $$r'(x,y)=\max \{r(x,y) , \min(r(x,a),r(b,y))\}.$$
In what follows,  we will  write $r'(x,y)=r(x,y)\vee(r(x,a)\wedge
r(b,y)).$ We claim that $r'$ is a Zadeh's fuzzy order on $X$ which
extends $r.$

\medskip

\noindent Claim 1. The fuzzy relation $r^{'}$  is Z--fuzzy
reflexive. Let $x \in X.$ Then, we have
$r'(x,x)=r(x,x)\vee(r(x,a)\wedge r(b,x))=1\vee(r(x,a)\wedge
r(b,x))=1.$ Thus, $r'$ is Z-fuzzy reflexive.

\medskip

\noindent Claim 2. The fuzzy relation $r^{'}$ is Z--fuzzy
antisymmetric. Indeed, let $x,y \in X$ such that $r'(x,y)>0$ and
$x\neq y.$ Then, as $r'(y,x)=r(y,x)\vee(r(y,a)\wedge r(b,x)),$ we
have four cases to study.

\medskip

\noindent First case. We have: $r'(x,y)= r(x,y)>0$ and
$r(x,a)\wedge r(b,y)>0.$ As $r$ is  Z--fuzzy antisymmetric,  we
get $r(y,x)=0.$ Hence, we obtain $r'(y,x)=r(y,a)\wedge r(b,x).$ On
the other hand since $r$ is Z--fuzzy transitive, so we get
$r(b,a)\geq r(b,y)\wedge r(y,a).$ Since $r(b,a)=0$ and $r(b,y)>0,$
then we obtain $r(y,a)=0.$ Thus, we have $r'(y,x)=0.$

\medskip

\noindent Second case. We have: $r'(x,y)= r(x,y)>0$ and
$r(x,a)\wedge r(b,y)=0.$ Then, as $r$ is Z--fuzzy antisymmetric,
we get $r(y,x)=0.$ Hence, $r'(y,x)=r(y,a)\wedge r(b,x).$ Since
$r(x,a)\wedge r(b,y)=0,$ so we have  $r(x,a)=0$ or $r(b,y)=0.$

\medskip

(a) First subcase. We have: $r(x,a)=0.$ Then, since $r$ is
Z--fuzzy transitive we get $r(x,a)\geq r(x,y)\wedge r(y,a).$ As
$r(x,a)=0$ and $r(x,y)>0,$ so we obtain $r(y,a)=0.$ Hence, we get
$r'(y,x)=0.$

\medskip

(b) Second subcase. We have: $r(b,y)=0.$ Then, as $r$ is Z--fuzzy
transitive,  we get $r(b,y)\geq r(b,x)\wedge r(x,y).$ Since
$r(b,y)=0$ and $r(x,y)>0,$ so we have $r(b,x)=0.$ Thus, we obtain
$r'(y,x)=0.$

\medskip

\noindent Third case. We have: $r'(x,y)= r(x,a)\wedge r(b,y)>0$
and $r(x,y)>0.$ Then,  as $r$ is Z--fuzzy antisymmetric we get
$r(y,x)=0.$ Hence, we obtain $r'(y,x)=r(y,a)\wedge r(b,x).$ Then,
since $r$ is Z--fuzzy transitive we get $r(b,a)\geq r(b,y)\wedge
r(y,a).$ As $r(b,a)=0$ and $r(b,y)>0,$ so we get $r(y,a)=0.$ Thus,
we obtain $r'(y,x)=0.$

\medskip

\noindent Fourth case. We have: $r'(x,y)= r(x,a)\wedge r(b,y)>0$
and $r(x,y)=0.$ So, as $r$ is Z--fuzzy transitive, we get
$r(b,a)\geq r(b,x)\wedge r(x,a).$ Since $r(b,a)=0$ and $r(x,a)>0,$
then we obtain $r(b,x)=0.$ Hence, we have $r(y,a)\wedge r(b,x)=0.$
So, we get $r'(y,x)=r(y,x).$ Then, by using  the Z--fuzzy
transitivity of $r$ we obtain $r(b,x)\geq r(b,y)\wedge r(y,x).$ As
$r(b,x)=0$ and $r(b,y)>0,$ so we obtain $r(y,x)=0.$ Hence, we get
$r'(y,x)=0.$ Therefore, $r'$ is Z--fuzzy antisymmetric.

\medskip

\noindent Claim 3. The fuzzy relation $r^{'}$ is Z--fuzzy
transitive. Indeed, let $x,y,z \in X$. Then, we have four cases to
study.

\medskip

\noindent First case. We have: $r^{'}(x,y)=r(x,y)$ and
$r^{'}(y,z)=r(y,z).$ Then we get $r^{'}(x,y)\wedge
r^{'}(y,z)=r(x,y)\wedge r(y,z).$ Since $r$ is Z-fuzzy transitive,
so we obtain $r(x,z)\geq r(x,y)\wedge r(y,z).$ Hence, $r(x,z)\geq
r^{'}(x,y)\wedge r^{'}(y,z).$ On the other hand since
$r'(x,z)=r(x,z)\vee(r(x,a)\wedge r(b,z)),$ then we get
$r'(x,z)\geq r(x,z).$ Thus, we have $r^{'}(x,z)\geq
r^{'}(x,y)\wedge r^{'}(y,z).$

\medskip

\noindent Second case. We have: $r^{'}(x,y)=r(x,a)\wedge r(b,y)$
and $r^{'}(y,z)=r(y,a)\wedge r(b,z).$ So we get, $r^{'}(x,y)\wedge
r^{'}(y,z)=r(x,a)\wedge r(b,y)\wedge r(y,a)\wedge r(b,z).$ Since
$r(b,a)=0$ and   $r$ is  Z-fuzzy transitive, so we get $r(b,a)\geq
r(b,y)\wedge r(y,a).$ Then, we obtain $r(b,y)\wedge r(y,a)=0.$
Hence, we have $r^{'}(x,y)\wedge r^{'}(y,z)=0.$ Thus,
$r^{'}(x,z)\geq r^{'}(x,y)\wedge r^{'}(y,z).$

\medskip

\noindent Third case. We have: $r^{'}(x,y)=r(x,y)$ and
$r^{'}(y,z)=r(y,a)\wedge r(b,z).$ Then, we get $r^{'}(x,y)\wedge
r^{'}(y,z)=r(x,y)\wedge r(y,a)\wedge r(b,z).$ On the other hand,
as $r'(x,z)=r(x,z)\vee(r(x,a)\wedge r(b,z)),$ so $r'(x,z)\geq
r(x,a)\wedge r(b,z).$ As $r$ is  Z-fuzzy transitive, then we get
$r(x,a)\geq r(x,y)\wedge r(y,a).$ Hence, we obtain $r'(x,z)\geq
r(x,y)\wedge r(y,a)\wedge r(b,z).$ Thus, we have $r^{'}(x,z)\geq
r^{'}(x,y)\wedge r^{'}(y,z).$

\bigskip

\noindent Fourth case. We have: $r^{'}(x,y)=r(x,a)\wedge r(b,y)$
and $r^{'}(y,z)=r(y,z).$ So we get, $r^{'}(x,y)\wedge
r^{'}(y,z)=r(x,a)\wedge r(b,y)\wedge r(y,z).$ On the other hand as
$r'(x,z)=r(x,z)\vee(r(x,a)\wedge r(b,z)),$ so we get $r'(x,z)\geq
r(x,a)\wedge r(b,z).$ In addition, By using the Z-fuzzy
transitivity of $r$ we get $r(b,z)\geq r(b,y)\wedge r(y,z).$
Hence, we obtain $r'(x,z)\geq r(x,a)\wedge r(b,y)\wedge r(y,z).$
So, we obtain $r^{'}(x,z)\geq r^{'}(x,y)\wedge r^{'}(y,z).$ Then,
we deduce that $r'(x,z) \geq \min\{r'(x,y),r'(y,z)\},$ for every
$y \in X.$ Thus, $r'$ is Z-fuzzy transitive. Therefore, $r'$ is a
Zadeh's fuzzy order on $X.$

\bigskip

\noindent Claim 4. The Zadeh's fuzzy order $r^{'}$ is an extension
of $r.$ Indeed, as for all $(x,y) \in X^{2},$ we have
$r'(x,y)=\max \{r(x,y), \min(r(x,a),r(b,y))\}\geq r(x,y).$ So,
$r'$ is an extension of $r.$  Moreover, we have
$$r'(a,b)=\max \{r(a,b), \min(r(a,a),r(b,b))\}=1$$
and
$$r'(b,a)=\max \{r(b,a),
\min(r(b,a),r(b,a))\}=r(b,a)=0.$$ \qed

\bigskip

Now we are able to give the proof of Theorem 2.1.

\bigskip

 \noindent{\bf Proof of Theorem 2.1.} Let $X=\{x_{1}, ...,
x_{n}\}$ be a finite set and let $r$ be a  Zadeh's fuzzy order on
$X$. If $x_{1}$ is comparable to every $x_{i}$ for $i\in \{2,...,
n\},$ it is ok. If note, by applying at most $n-1$ times Lemma
2.2, we get a Zadeh's fuzzy order $r_{1},$ say, which  extends $r$
and satisfies
$$r_{1}(x_{1}, x_{i}) > 0 \hbox{ or } r_{1}(x_{i}, x_{1}) > 0 \hbox{ for every }
i \in \{1, 2, ..., n\}.$$ Now, if  $x_{2}$ is comparable to every
$x_{i}$ for $i\in \{3,..., n\},$ it is ok. If it is note the case,
by applying  at most $n-2$ times  Lemma 2.2, we get a Zadeh's
fuzzy order $r_{2}$ which  extends $r_{2}$ such that
$$r_{2}(x_{j}, x_{i}) > 0 \hbox{ or } r_{2}(x_{i}, x_{j}) > 0 \hbox{ for every }
i \in \{1, 2,..., n\}\hbox{ and }j=1, 2.$$
 By induction for every $k \in \{1,..., n-1\},$ there exists a Zadeh's fuzzy
order $r_{k},$ say, which  extends $r_{k-1}$  and satisfies the
following:
$$r_{k}(x_{j}, x_{i}) > 0 \hbox{ or }r_{k}(x_{i}, x_{j}) > 0 \hbox{ for every }
j \in \{1, 2, ..., n\}\hbox{ and }i=1, 2, ..., k.$$ Since $r\leq
r_{1} \leq  r_{2} \leq ...\leq r_{n-2} \leq r_{n-1},$ then $r \leq
r_{n-1}.$ Thus, we have obtain a linear Zadeh's fuzzy order
$r_{n-1}$ which extends $r.$ \qed

\begin{rem} Let $k$ be number of the times of application of
Lemma 2.2 in the proof of Theorem 2.1. Then, the naturel number
$k$ satisfies $k\leq \frac{m}{2}$ where $m=card(A)$ and
$$A=\{(x_{i},x_{j})\in X^{2},~~\mbox {such that}~~
r(x_{i},x_{j})=r(x_{j},x_{i})=0\}.$$ As $m \leq n(n-1),$  then we
get $k \leq \frac{n(n-1)}{2}.$ \end{rem}

\section{A characterization of Zadeh's fuzzy orders by their linear
extensions} In this section, we will characterize  Zadeh's fuzzy
orders which are defined on a finite nonempty set $X$ by their
linear extensions. More precisely, we will prove the following
result.

\begin{thm}
Every Zadeh's fuzzy order on a finite nonempty set $X$ is the
fuzzy intersection of all linear Zadeh's fuzzy orders which extend
it.
\end{thm}

In order to prove Theorem 3.1, we will need the following
technical lemma.

\begin{lem}  Let $(X, r)$ be a finite nonempty Zadeh's
 fuzzy ordered set and let $a, b$ be two elements of $X$ such that
 $r(a, b) > 0.$ Then, there exists at least a linear Zadeh's fuzzy order $r'$ on $X$
which extends $r$ such that $r'(a,b)=r(a, b).$
\end{lem}

\noindent{\bf Proof}. Let $(X, r)$ be a nonempty Zadeh's fuzzy
ordered set and let $a, b$ be two elements of $X$ such that $r(a,
b) > 0.$ If $r$ is linear then we take $s=r$ and Lemma 3.2 is
proved. If it is note the case, then from Theorem 2.1 there exist
a linear Zadeh's fuzzy order $r^{'}$ on $X$ which extends $r.$ If
$r(a,b)= r^{'}(a, b),$ then we take $s=r^{'}$ and  Lemma 3.2 is
proved. If it note the case, then assume that we have $r(a,b)<
r^{'}(a, b).$ So, we set $\beta=r(a,b).$ Hence, we get $0 < \beta
< 1.$ Next, we will define the following fuzzy relation $s$ on $X$
by setting:

$$s(x,y)=\left\{\begin{array}{c}
r^{'}(x, y), ~~~~~if  ~~ r(x,y)> \beta ; \\
\beta \wedge r^{'}(x, y), ~if  ~~ r(x,y)\leq \beta.\\
\end{array}\right.$$

 \noindent We claim that $s$ is a Zadeh's fuzzy order which extends $r$ and
satisfies $s(a, b)=r(a, b).$

\bigskip

\noindent Claim 1. The fuzzy relation $s$ is a linear Zadeh's
fuzzy order on $X$.

\medskip

\noindent (i) Z--fuzzy reflexivity. For all $x \in X$ we have
$r(x,x)=1>\beta$, then we get $s(x,x)=r^{'}(x,x)=1.$ Thus $s$ is
fuzzy  reflexive.

\medskip

\noindent (ii) Z--fuzzy antisymmetry. Let $x,y \in X$ such that
$s(x,y)>0$ and $x\neq y.$ To show that $s(y, x)=0,$ we have to
consider two subcases.

\medskip

First subcase. We have: $s(x,y)=r^{'}(x,y)>0.$ Then,  $r(x,y)>
\beta.$ As $r^{'}$ is Z--fuzzy antisymmetric, so we get
$r^{'}(y,x)=0$. On the other hand, since $r(x,y)> \beta>0$ and $r$
is antisymmetric, then we get $r(y,x)=0\leq \beta.$ Hence,
$s(y,x)=\beta \wedge r^{'}(y, x)=0$. Thus, $s$ is fuzzy
antisymmetric.

\medskip

Second subcase. We have: $s(x,y)=\beta \wedge r^{'}(x, y)>0.$
Then,  we get  $r^{'}(x, y)>0.$ By the fuzzy
 antisymmetry of $r^{'}$ we obtain $r^{'}(y,x)=0.$ As $r(y, x) \leq
 r^{'}(y, x),$ then we obtain $r(y, x)=0 \leq \beta.$ Hence, we get
  $s(y,x)=\beta \wedge r^{'}(y,x)=0.$ Thus, $s$ is fuzzy  antisymmetric.

\medskip

\noindent (iii) Z--fuzzy transitivity. Let $x,y,z \in X$. To prove
that $s(x, z) \geq \{s(x, y), s(y, z) \},$ we have to distinguish
the following four cases.

\medskip

\noindent First case. We have: $s(x,y)=r^{'}(x,y)$ and
$s(y,z)=r^{'}(y,z).$ Then, we get $r(x,y)> \beta$ and $r(y,z)>
\beta$. As $r$ is  Z--fuzzy transitive, so $r(x,z)\geq
r(x,y)\wedge r(y,z).$ Hence, $r(x,z)> \beta$. Then, we obtain
$s(x,z)=r^{'}(x,z)$. From the fuzzy transitivity of $r^{'}$ we get
$r^{'}(x,z)\geq r^{'}(x,y)\wedge r^{'}(y,z).$ Thus, $s(x,z)\geq
s(x,y)\wedge s(y,z).$

\medskip

\noindent Second case. We have: $s(x,y)=r^{'}(x,y)$ and
$s(y,z)=\beta \wedge r^{'}(y,z).$ Then, we get  $r(x,y)> \beta$
and $r(y,z)\leq \beta.$ So, we obtain $r(x,y)\wedge r(y,z)\leq
\beta.$ Next, we have two subcases to consider.

\medskip

First subcase. We have: $r(x,z)> \beta.$ So, we get
$s(x,z)=r^{'}(x,z).$ As $r^{'}$ is Z--fuzzy transitive,
 we obtain $r^{'}(x,z)\geq r^{'}(x,y)\wedge r^{'}(y,z).$ So, we
 get $r^{'}(x,z)\geq r^{'}(x,y)\wedge (\beta \wedge r^{'}(y,z)).$ Thus,
we obtain $s(x,z)\geq s(x,y)\wedge s(y,z).$

\medskip

Second subcase. We have: $r(x,z)\leq \beta.$ Then, $s(x,z)=\beta
\wedge r^{'}(x,z)$. Since  $r^{'}$ is  Z--fuzzy transitive, hence
 we get $r^{'}(x,z)\geq r^{'}(x,y)\wedge r^{'}(y,z)$ and
 $\beta \wedge r^{'}(x,z)\geq r^{'}(x,y)\wedge
(\beta \wedge r^{'}(y,z)).$ Thus, $s(x,z)\geq s(x,y)\wedge
s(y,z)$.

\medskip

Third case. We have: $s(x,y)=\beta \wedge r^{'}(x,y)$ and
$s(y,z)=r^{'}(y,z).$ Then, we get  $r(x,y)\leq \beta$ and $r(y,z)>
\beta.$ So, $r(x,y)\wedge r(y,z)\leq \beta.$ For this case we have
two subcases to study.

\medskip

First subcase. We have: $r(x,z)> \beta.$ Then, we get
$s(x,z)=r^{'}(x,z).$ Since $r^{'}$ is Z--fuzzy transitive, so we
have $r^{'}(x,z)\geq r^{'}(x,y)\wedge r^{'}(y,z)$ and
$r^{'}(x,z)\geq (\beta \wedge r^{'}(x,y))\wedge r^{'}(y,z).$ Thus,
we obtain $s(x,z)\geq s(x,y)\wedge s(y,z).$

\medskip

Second subcase. We have: $r(x,z)\leq \beta.$ Then, $s(x,z)=\beta
\wedge r^{'}(x,z).$ Since $r^{'}$ is Z--fuzzy transitive, we get
$r^{'}(x,z)\geq r^{'}(x,y)\wedge r^{'}(y,z)$ and $\beta \wedge
r^{'}(x,z)\geq (\beta \wedge r^{'}(x,y))\wedge r^{'}(y,z).$ Thus,
$s(x,z)\geq s(x,y)\wedge s(y,z).$

\medskip

Fourth case. We have: $s(x,y)=\beta \wedge r^{'}(x,y)$ and
$s(y,z)=\beta \wedge r^{'}(y,z).$ Then, we get  $r(x,y)\leq \beta$
and $r(y,z)\leq \beta.$ So, $r(x,y)\wedge r(y,z)\leq \beta$. For
this case we have  two subcases to distinguish.

\medskip

First subcase. We have: $r(x,z)> \beta.$ Then, we have
$s(x,z)=r^{'}(x,z).$ As $r^{'}$ is fuzzy transitive, so we obtain
$r^{'}(x,z)\geq r^{'}(x,y)\wedge r^{'}(y,z).$ Hence,
$r^{'}(x,z)\geq (\beta \wedge r^{'}(x,y))\wedge (\beta \wedge
r^{'}(y,z)).$ Thus, we get $s(x,z)\geq s(x,y)\wedge s(y,z).$

\medskip

Second subcase. We have: $r(x,z)\leq \beta.$ Then, we get
$s(x,z)=\beta \wedge r^{'}(x,z).$ Since $r^{'}$ is
 Z--fuzzy transitive, so we obtain $r^{'}(x,z)\geq r^{'}(x,y)\wedge
 r^{'}(y,z)$ and $\beta \wedge r^{'}(x,z)\geq (\beta \wedge r^{'}(x,y))\wedge (\beta \wedge r^{'}(y,z)).$
 Thus, we get $s(x,z)\geq s(x,y)\wedge s(y,z).$ Thus, we obtain
 $s(x,z) \geq \min\{s(x,y),s(y,z)\}$,  for all $y \in X.$ Hence,
 $s$ is fuzzy transitive. Therefore, $s$ is a Zadeh's
fuzzy order on $X.$

\medskip

\noindent Claim 2. The Z--fuzzy order relation $s$ is linear.
Indeed, let $x, y \in  X,$ such that $x\neq y$. Then, since
$r^{'}$ is a linear Zadeh's fuzzy order we get  $r^{'}(x,y)>0$ or
$r^{'}(y,x)>0.$

First case. We have: $r^{'}(x,y)>0.$ Then, $s(x, y)=r^{'}(x,y)$ or
$s(x, y)=\beta \wedge r^{'}(x, y).$ So, we obtain $s(x, y) > 0.$

Second case. We have: $r^{'}(y, x)>0.$ Then, $s(y, x)=r^{'}(y, x)$
or $s(y, x)=\beta \wedge r^{'}(y,x).$ So, we get $s(y, x) > 0.$

\bigskip

 \noindent Claim 3. We have: $r(x,y)\leq s(x,y).$ To show this,
we will distinguish two cases.

First case. We have: $r(x,y)> \beta.$  Then, $s(x,y)=r^{'}(x, y).$
As $r^{'}$ is an extension of $r,$ so $r(x,y) \leq s(x, y).$

Second case. We have: $r(x,y)\leq \beta.$ Then, $r^{'}(x, y)=\beta
\wedge r^{'}(x, y).$ Since $r(x,y)\leq \beta$ and $r(x, y) \leq
r^{'},$ then $r(x, y) \leq \beta \wedge r^{'}(x, y).$ Thus, we
have $r(x, y) \leq s(x,y).$ Therefore, the Lemma 3.2 is
proved.\qed

\bigskip

Next, we are ready to give the proof of Theorem 3.1.

 \noindent{\bf Proof of Theorem 3.1.} Let $r$ be a Zadeh's fuzzy order on a
nonempty finite set $X.$ If $r$ is a linear Zadeh's fuzzy order,
then it is ok. If note, let $A$ be the set of all linear Zadeh's
fuzzy orders on $X$ which extend $r$. From Theorem 2.1,
 $A$ is nonempty. Let $s_{0}$ be the fuzzy intersection of all elements of $A.$
Then, we have $$s_{0}(x, y)=\inf_{s \in A} s(x, y),\hbox{ for
every }(x, y) \in X^{2}.$$

\bigskip

\noindent Claim 1. The fuzzy relation $s_{0}$ is a Zadeh's fuzzy
order on $X.$

(a) Z--fuzzy reflexivity. Let $x \in X.$ Then, $s(x, x)=1$ for
every $s \in A.$ So, $s_{0}(x, x)=1.$ Thus, $r$ is Z--fuzzy
reflexive.

(b) Z--fuzzy antisymmetry. Let $x, y \in X$ such that $x\not=y$
and $s_{0}(x, y) > 0.$ Let $s \in A.$ Then, $s(x, y) > 0.$ On the
other hand, we know that $s$ is Z--fuzzy antisymmetric. So, we get
$s(y, x)=0.$ Hence, $s_{0}(y, x)=\inf_{s \in A} s(y, x)=0.$ Thus,
$s_{0}$ is Z--fuzzy antisymmetric.

(c) Z--fuzzy transitivity. Let $x, y, z \in X$ and $s\in A.$ Then,
we have
$$s_{0}(x, y) \leq s(x, y) \hbox{ and } s_{0}(y, z) \leq s(y,
z).$$ Hence, we get
$$\min\{s_{0}(x, y), s_{0}(y, z)\} \leq \min\{s(x, y), s(y,
z)\}.$$ On the other hand, we know that $s$ is Z--fuzzy
transitive. So, we get $$\min\{s(x, y), s(y, z)\} \leq s(x, z).$$
Then, we obtain, $$\min\{s_{0}(x, y), s_{0}(y, z)\} \leq s(x, z),
\hbox{ for every } s \in A.$$ Since $s_{0}(x, z)=\inf_{s \in A}
s(x, z),$ so we get $$\min\{s_{0}(x, y), s_{0}(y, z)\} \leq
s_{0}(x, z).$$ Thus, $s_{0}$ is Z--fuzzy transitive. Therefore,
$s_{0}$ is a Zadeh's fuzzy order on $X.$

\bigskip

\noindent Claim 2. We have: $r \leq s_{0}.$ Indeed if $x, y \in
X,$ then we get
$$r(x, y) \leq s(x, y),\hbox{ for every } s \in A.$$
So, we obtain
$$r(x, y) \leq s_{0}(x, y),\hbox{ for every } (x, y)\in X^{2}.$$
Thus, we have $r \leq s_{0}.$

\bigskip

\noindent Claim 3. For every $a, b \in X,$ we have
$$( r(a, b)=r(b, a)=0 ) \Rightarrow ( s_{0}(a, b)=s_{0}(b, a)=0
).$$ Indeed, assume that $r(a, b)=r(b, a)=0.$ Then, from Theorem
2.1, there exists two linear Zadeh's fuzzy orders of $r,$ $s_{1}$
and $s_{2},$ say such that
$$ \left\{ \begin{array}{c}
s_{1}(a, b)=1  \\
s_{1}(b, a)=0
\end{array}
\right.$$ and
$$\left\{\begin{array}{c}
s_{2}(a, b)=0  \\
s_{2}(b, a)=1.
\end{array}
\right. $$ \\ As $s_{0} \leq s_{1}$ and $s_{0} \leq s_{2},$ hence
we get $s_{0}(a, b)=s_{0}(b, a)=0.$

\noindent Claim 4. We have: $r=s_{0}.$ Indeed, we know that $r\leq
s_{0}.$ By absurd assume that $r\not=s_{0}.$ Then, there is $a, b
\in X$ such that $r(a, b) < s_{0}(a, b).$ Hence, we get $s_{0}(a,
b)> 0.$ As $s_{0}$ is Z-fuzzy antisymmetric, so we obtain
$s_{0}(b, a)=0.$ On the other hand,  we know that $r \leq s_{0}.$
Then, $r(b, a)=0.$ Now, we claim that $r(a, b) > 0.$ On the
contrary assume that $r(a, b)=0.$ So, we get $r(a, b)=r(b, a)=0.$
By using Claim 2, we obtain $s_{0}(a, b)=s_{0}(b, a)=0.$ That is a
contradiction with the fact that $s_{0}(a, b) > 0.$ As $r(a, b)>
0,$ then from Lemma 3.2, there exists $s \in A$ such that $s(a,
b)=r(a, b).$ Since $s_{0} \leq s,$ hence we get $s_{0}(a, b) \leq
r(a, b).$  That is a contradiction with our assumption that $r(a,
b) < s_{0}(a, b).$ Therefore, we obtain $r=s_{0}.$ \qed

\section{Examples of construction of linear extensions of Zadeh's fuzzy orders }
In this section we will give some examples of construction of
linear extensions of Zadeh's fuzzy orders by applying Theorem 2.1.

\noindent 1. Let $X=\{a, b, c\}$ and $r$ be the Zadeh's fuzzy
order on $X$ defined by the following table

\begin{center}

\begin{tabular}{|c|c|c|c|}
  \hline
& a & b & c  \\
 \hline
a & 1 & 0 & $\gamma$ \\
b & 0 & 1 & 0 \\
c & 0 & 0 & 1 \\
  \hline
\end{tabular}
\end{center}

\noindent By using Lemma 2.2, we get in the first step a Zadeh's
fuzzy order $r_{1}$  which extends $r$ and it is  defined by
$r_{1}(x,y)=\max \{r(x,y) , \min(r(x,a),r(b,y))\}$.  Then, we
represent $r_{1}$ by the following table

\begin{center}

\begin{tabular}{|c|c|c|c|}
  \hline
& a &b & c \\
 \hline
a & 1 & 1 & $\gamma$ \\
b & 0 & 1 & 0 \\
c & 0 & 0 & 1 \\
  \hline
\end{tabular}
\end{center}

Since $r_{1}$ is not linear, then we apply again Lemma 2.2  for
$r_{1}.$ So, we get a Zadeh's fuzzy order $r_{2}$  extending
$r_{1}$ (also $r_{1}$ extends also $r$) defined by
$$r_{2}(x,y)=\max \{r_{1}(x,y) , \min(r_{1}(x,b),r_{1}(c,y))\}.$$ We
can represent $r_{2}$ by the following table

\begin{center}

\begin{tabular}{|c|c|c|c|}
  \hline
& a &b & c \\
 \hline
a & 1 & 1 & 1 \\
b & 0 & 1 & 1 \\
c & 0 & 0 & 1 \\
  \hline
\end{tabular}
\end{center}

The procedure stop here because $r_{2}$ is a linear Zadeh's fuzzy
order extending $r$.

In this example the number of application of Lemma 2.2 is $k=2.$
In this case, we have $k=\frac{m}{2}$ such that $m=card(A)$ with
$A=\{(a,b), (b,a), (b,c), (c,b)\}.$

\medskip

\noindent 2. Let $X=\{a, b, c, d\}$ and $r$ be the Zadeh's fuzzy
order on $X$ defined by the following matrix:

\begin{center}

\begin{tabular}{|c|c|c|c|c|}
  \hline
& a & b & c & d \\
 \hline
a & 1 & 0 & $\gamma$ & $\lambda$ \\
b & 0 & 1 & $\alpha$ & $\beta$ \\
c & 0 & 0 & 1 & 0 \\
d & 0 & 0 & 0 & 1 \\
  \hline
\end{tabular}
\end{center}

\noindent By using Lemma 2.2, we get in the first step a Zadeh's
fuzzy order $r_{1}$ which extends $r$ and it is defined by
$r_{1}(x,y)=\max \{r(x,y) , \min(r(x,a),r(b,y))\}.$  Hence, we
represent $r_{1}$ by the following table:

\begin{center}

\begin{tabular}{|c|c|c|c|c|}
  \hline
& a & b & c & d \\
 \hline
a & 1 & 1 & $\gamma \vee \alpha$ & $\lambda \vee \beta$ \\
b & 0 & 1 & $\alpha$ & $\beta$ \\
c & 0 & 0 & 1 & 0 \\
d & 0 & 0 & 0 & 1 \\
  \hline
\end{tabular}
\end{center}

As $r_{1}$ is not linear, then we can still apply Lemma 2.2 for
$r_{1}.$ So, we obtain a  Zadeh's fuzzy order $r_{2},$ say, which
extends $r_{1}$ and it is  defined by
$$r_{2}(x,y)=\max \{r_{1}(x,y) , \min(r_{1}(x,c),r_{1}(d,y))\}.$$ Then, we
represent $r_{2}$ by the following table:

\begin{center}

\begin{tabular}{|c|c|c|c|c|}
  \hline
& a & b & c & d \\
 \hline
a & 1 & 1 & $\gamma \vee \alpha$ & $\lambda \vee \beta \vee \gamma \vee \alpha$ \\
b & 0 & 1 & $\alpha$ & $\beta \vee \alpha$ \\
c & 0 & 0 & 1 & 1 \\
d & 0 & 0 & 0 & 1 \\
  \hline
\end{tabular}
\end{center}

The procedure stop here because $r_{2}$ is a linear Zadeh's fuzzy
order extending $r$. In this example the number of application of
Lemma 2.2 is $k=2.$ Thus,  we have $k=\frac{m}{2}$ where
$m=card(A)$ with $A=\{(a,b), (b,a), (c,d), (d,c)\}.$

\bigskip

\noindent 3. Let $X=\{x_{1}, x_{2}, x_{3}, x_{4}, x_{5}, x_{6},
x_{7}\}$ and $r$ be the Zadeh's fuzzy order on $X$ defined by the
following matrix:

\begin{center}

\begin{tabular}{|c|c|c|c|c|c|c|c|}
  \hline
& $x_{1}$ &$ x_{2}$ & $x_{3}$ & $x_{4}$ & $x_{5}$& $x_{6}$& $x_{7}$\\
 \hline
$x_{1}$ & 1 & 0 & 0 & 0.55 & 0.40 &0.45 &0.60\\
$x_{2}$ & 0 & 1 & 0 & 0.60 & 0.50 & 0.35 &0.75\\
$x_{3}$ & 0.15 & 0 & 1 & 0.30 & 0.70 &0.80 &0.90\\
$x_{4}$ & 0 & 0 & 0 & 1 & 0 &0.15 &0\\
$x_{5}$ & 0 & 0 & 0 & 0 & 1 &0.30 &0.25\\
$x_{6}$ & 0 & 0& 0 & 0 & 0 &1 &0\\
$x_{7}$ & 0 & 0& 0 & 0 & 0 &0.20 &1\\
  \hline
\end{tabular}
\end{center}

\noindent By applying a finite number of times Lemma 2.2, we
obtain a linear extension of $r$ which we represented by the
following matrix

\begin{center}

\begin{tabular}{|c|c|c|c|c|c|c|c|}
  \hline
& $x_{1}$ &$ x_{2}$ & $x_{3}$ & $x_{4}$ & $x_{5}$& $x_{6}$& $x_{7}$\\
 \hline
$x_{1}$ & 1 & 1 & 0 & 0.60 & 0.60 &0.45 &0.75\\
$x_{2}$ & 0 & 1 & 0 & 0.60 & 0.60 & 0.35 &0.75\\
$x_{3}$ & 0.15 & 0.15 & 1 & 0.30 & 0.70 &0.80 &0.90\\
$x_{4}$ & 0 & 0 & 0 & 1 & 1 &0.30 &0.25\\
$x_{5}$ & 0 & 0 & 0 & 0 & 1 &0.30 &0.25\\
$x_{6}$ & 0 & 0& 0 & 0 & 0 &1 &0\\
$x_{7}$ & 0 & 0& 0 & 0 & 0 &0.20 &1\\
  \hline
\end{tabular}
\end{center}

\end{document}